\input amstex
\documentstyle{amsppt}
\magnification=\magstep1
\pagewidth{6.5truein}
\pageheight{9.0truein}
\NoBlackBoxes
\def\ep{\varepsilon}
\def\D{{\Cal D}}
\def\L{{\Cal L}}
\def\N{{\Cal N}}
\def\U{{\Cal U}}
\def\bs{\text{\bf s}}
\def\nat{{\Bbb N}}
\def\real{{\Bbb R}}

\def\un#1{\underline{#1}}
\def\eqdf{\mathop{\buildrel {\text{df}}\over =}\nolimits}
\def\diam{\operatorname{diam}}
\def\dist{\operatorname{dist}}
\topmatter
\title On wide-$(s)$ sequences and their applications to certain
classes of operators\endtitle
\rightheadtext{On WIDE-$(s)$ SEQUENCES}
\author H. Rosenthal\endauthor
\affil Department of Mathematics\\ The University of Texas at Austin\\
Austin, TX 78712-1082\endaffil
\date August 5, 1996\enddate
\thanks Research partially supported by NSF Grant DMS-9500874 and
TARP Grant ARP-275.\endthanks
\abstract
A basic sequence in a Banach space is called wide-$(s)$ if it is bounded and
dominates the summing basis. (Wide-$(s)$ sequences were originally
introduced by I.~Singer, who termed them $P^*$-sequences).
These sequences and their quantified versions,
termed $\lambda$-wide-$(s)$ sequences, are used to characterize various
classes of operators between Banach spaces, such as the weakly compact,
Tauberian, and super-Tauberian operators,
as well as a new intermediate class introduced here,
the strongly Tauberian operators.
This is a nonlocalizable class which nevertheless forms an open semigroup and
is closed under natural operations such as taking double adjoints.
It is proved for example that an operator is non-weakly compact iff for every
$\ep >0$, it maps some $(1+\ep)$-wide-$(s)$-sequence to a wide-$(s)$ sequence.
This yields the quantitative triangular arrays result
characterizing reflexivity, due to R.C.~James.
It is shown that an operator is non-Tauberian (resp. non-strongly Tauberian)
iff for every $\ep>0$, it maps some  $(1+\ep)$-wide-$(s)$ sequence into
a norm-convergent sequence (resp. a sequence whose image has diameter less
than $\ep$). This is applied to obtain a direct ``finite'' characterization of
super-Tauberian operators, as well as the following characterization, which
strengthens a recent result of M.~Gonz\'alez and A.~Mart{\'\i}nez-Abej\'on:
An operator is non-super-Tauberian iff there are for every $\ep>0$,
finite $(1+\ep)$-wide-$(s)$ sequences of arbitrary length whose images have
norm at most $\ep$.
\endabstract
\endtopmatter

\document
\baselineskip=18pt		

\subhead \S1. Introduction\endsubhead

A semi-normalized basic sequence $(b_j)$ in a Banach space is called
wide-$(s)$ if it dominates the summing basis; i.e., $\sum c_j$ converges
whenever $\sum c_jb_j$ converges.
Our main objective here is to show that wide-$(s)$ sequences provide a
unified approach for dealing with certain important classes of operators
between Banach spaces.
For example, we show in Proposition~2 that an operator between Banach spaces
is non-weakly compact iff it maps some wide-$(s)$ sequence into a
wide-$(s)$ sequence.
Similarly, we show in Theorem~5 that an operator is non-Tauberian iff it
maps some wide-$(s)$ sequence into a norm convergent sequence.
In Corollary~6 we obtain that an operator is Tauberian iff it maps some
subsequence of a given wide-$(s)$ sequence, into a wide-$(s)$ sequence.
(Recall that $T:X\to Y$ is Tauberian if $T^{**} (X^{**}\sim X) \subset 
Y^{**} \sim Y$.)
Theorem~5 may also be deduced from the results due to R.~Neidinger and the
author in \cite{NR}, and standard facts.
However we give here a self-contained treatment.

After circulating the first version of this paper, we learned that
wide-$(s)$ sequences were originally introduced by I.~Singer, who termed
them ``$P^*$-sequences'' \cite{S1}.
Some of the pleasant permanence properties discovered by Singer, are
as follows:
a semi-normalized sequence $(x_j)$ is wide-$(s)$ iff its difference
sequence $(x_{j+1}-x_j)$ is a basic sequence.
Moreover if $(x_j)$ is a semi-normalized basis for a Banach space $X$ and
$(x_j^*)$ is its sequence of biorthogonal functionals, then $(x_j)$ is
wide-$(s)$ iff $(\sum_{j=1}^n x_j^*)_{n=1}^\infty$ is a wide-$(s)$
sequence in $X^*$ (cf.\ Theorem~9.2, page 311 of \cite{S2}; as noted in
our earlier version, these equivalences also follow from certain arguments in
\cite{R3}).
Of course convex block bases of wide-$(s)$ sequences are also wide-$(s)$
(see Proposition~3 below).

We prove in Corollary~17 that suitable perturbations of wide-$(s)$
sequences have wide-$s$ subsequences.
In fact our argument yields that given $(x_j)$ a wide-$(s)$ sequence in
a Banach space $X$, there exists an $\ep>0$ so that if $(y_j)$ is any
bounded perturbation of $(x_j)$ so that all $w^*$-cluster points of
$(x_j-y_j)$ in $X^{**}$ have norm at most $\ep$, then $(y_j)$ has a
wide-$(s)$ subsequence.
We also obtain that $\ep$ may be chosen depending only on $\lambda$ the
wide-$(s)$ constant of $(x_j)$, which we  introduce in Definition~3.
It follows immediately from this definition that if $(b_j)$ is a
$\lambda$-wide-$(s)$ sequence in $X$, then $\|b_j\| \le\lambda$ for all $j$
and there exists a sequence $(f_j)$ in $X^*$ with $\|f_j\|\le \lambda$
for all $j$ so that
$$\align
f_i(b_j) & = 1\ \text{ all }\ 1\le i\le j<\infty \\
f_i(b_j) & =0\ \text{ all }\ 1\le j<i<\infty\ .
\endalign$$
(We call sequences $(b_j)$ satisfying this condition  {\it triangular\/},
because the matrix $(f_i(b_j))$ is obviously upper triangular consisting
of 1's above and on the natural diagonal, zeros below.)
In Theorem~11 we prove that every non-reflexive Banach space has a
$(1+\ep)$-wide-$(s)$ sequence for every $\ep>0$; this immediately yields
the remarkable quantitative information on triangular arrays in
non-reflexive spaces discovered
by R.C.~James \cite{J1}, \cite{J2} and D.P.~Milman--V.D. Milman \cite{MM}
(cf.\ the remark following the statement of Theorem~11).

Our proof of this result involves the rather technical result
Theorem~12, and some standard facts on basic-sequence selections, formulated
in \cite{R3} and repeated here for completeness as Lemmas~13 and 14.
These facts are used to deduce the rather surprising result that given
$\ep>0$ and any  bounded non-relatively compact sequence $(x_j)$ in a
Banach space $X$, there is a subsequence $(x'_j)$ of $(x_j)$ and an $x$ in $X$
with $(b_j)$ a $(2+\ep)$-basic sequence, where $b_j=x'_j-x$ for all $j$.
Moreover, if $(x_j)$ is non-relatively weakly compact, there is a $c>0$ 
and a convex block basis $(u_j)$ of $(b_j)$ so that $(cu_j)$ is
$(1+\ep)$-wide-$(s)$ (Corollary~15).
(For an interesting recent result on uniformity in the biorthogonal constant
of subsequence of uniformly separated sequences see \cite{HKPTZ}.) 
Corollary~15 yields also a quantitative refinement of both Theorem~11
and Proposition~2a: given $\ep>0$, a non-weakly compact operator maps some
$(1+\ep)$-wide-$(s)$ sequence into a wide-$(s)$ sequence
(Theorem~11$'$, given in the second remark following the proof of
Corollary~15).

Our definition of $\lambda$-wide-$(s)$ sequences applies to finite sequences
as well.
This leads to localizations of our results.
For example, Proposition~19 yields that an operator is non-super weakly
compact iff there is a $\lambda\ge1$ so that for all $n$, it maps some
$\lambda$-wide-$(s)$ sequence of length $n$ into a $\lambda$-wide-$(s)$
sequence.
Applying a quantitative refinement of Theorem~5 (called Theorem~$5'$ and
formulated in the first remark following the proof of Corollary~15), we
obtain that a Banach space $X$ is non-super reflexive iff it contains
$(1+\ep)$-wide-$(s)$ sequences of arbitrarily large length, for every $\ep>0$.
In Proposition~20, we localize Theorem~5 to the setting of super Tauberian
operators.
Our results here are motivated by recent work of M.~Gonz\'alez and
A.~Mart\'\i nez-Abej\'on \cite{GM}.
In fact, Proposition~20 may be deduced from results in \cite{GM} and our
localization result, Proposition~18, whose proof does not require
ultraproduct techniques.
However using ultraproducts, we obtain certain of the results in \cite{GM}
by localizing our Theorem~$5'$ and Proposition~3.
Thus we obtain an alternative route to the result obtained in \cite{GM}:
if $T$ is a given operator, then if $T$ is non-super Tauberian,
and $T_U$ is an ultrapower of $T$, then $\ker T_U$ is non-reflexive.
We then obtain the following strengthening of Proposition~12 of \cite{GM}:
{\it An operator $T$ between Banach spaces is non-super-Tauberian
iff for all $\ep>0$, all $n$, $T$ maps some $(1+\ep)$-wide-$(s)$
sequence of length $n$ into a sequence whose elements have norm at
most $\ep$\/}.
(See the remark following the proof of Proposition~20.)
Our Proposition~20 also yields the immediate Corollary that the
super-Tauberian operators in  $\L(X,Y)$ are an open set.
(This result is due to D.G.~Tacon \cite{T}; the deduction of this  result
in \cite{GM} motivates our formulation of Proposition~20.)

Next we introduce a class of operators intermediate between Tauberian
and super Tauberian; the strongly Tauberian operators.
These are operators $T\in \L(X,Y)$ whose natural induced map
$\tilde T:X^{**}/X \to Y^{**}/Y$ is an isomorphism.
It is immediate that the strongly Tauberian operators are an open subset
of $\L(X,Y)$ and of course have the semigroup property.
We give several equivalences in Theorem~21, obtaining for example, that
$T$ is strongly Tauberian iff $T^{**}$ is strongly Tauberian iff $T$ maps
some subsequence of a given $\lambda$-wide-$(s)$ sequence into a
$\beta$-wide-$(s)$ sequence, where $\beta$ depends only on $\lambda$.
The proof also yields that $T$ is non-strongly Tauberian iff given
$\ep>0$, $T$ maps some $(1+\ep)$-wide-$(s)$ sequence into a sequence of
diameter less than $\ep$ (Proposition~22).
Of course this yields immediately, via Proposition~20, that every super
Tauberian operator is strongly Tauberian.
We then use these discoveries about strongly Tauberian operators in an
essential way, to localize Corollary~6.
Thus we obtain in Corollary~25 that $T$ is super Tauberian iff for all $k$
every $\lambda$-wide-$(s)$ sequence of length $n$ has a subsequence of
length $k$ mapped to a $\beta$-wide-$(s)$ sequence by $T$, where $\beta$
depends only on $\lambda$, $n$ only on $\lambda$ and $k$.

We conclude with a rather delicate localization of the infinite
perturbation result given in Corollary~17.
This result, Proposition~26, apparently requires ultraproduct techniques
for its proof.
For standard facts about ultraproducts in Banach spaces, which we use
without explicit reference, see \cite{H}.

\subhead \S2\endsubhead

We first summarize some of the basic concepts used here.

\definition{Definition 1}
A semi-normalized sequence $(b_j)$ in a Banach space is called
\roster
\item"(i)" a wide-$(s)$ sequence if $(b_j)$ is a basic sequence which
dominates the summing basis; i.e., $\sum c_j$ converges whenever
$\sum c_jb_j$ converges.
\item"(ii)" an $(s)$-sequence if $(b_j)$ is weak-Cauchy and a
wide-$(s)$ sequence.
\item"(iii)" an $\ell^1$-sequence if $(b_j)$ is equivalent to the usual
$\ell^1$-basis.
\item"(iv)" non-trivial weak-Cauchy if $(b_j)$ is weak-Cauchy but not
weakly convergent.
\endroster
\enddefinition

As proved in Proposition~2.2 of \cite{R3}, every non-trivial
weak-Cauchy sequence has an $(s)$-subsequence.
It then follows immediately from the $\ell^1$-theorem that

\proclaim{Proposition 1}
Every wide-$(s)$ sequence has a subsequence which is either an
$(s)$-sequence or an $\ell^1$-sequence.
\endproclaim

\noindent
(The $\ell^1$-theorem refers to the author's result that every bounded
non-trivial weak-Cauchy sequence has an $\ell^1$-subsequence.)

Of course every sequence which is either an $(s)$-sequence, or an
$\ell^1$-sequence, is wide-$(s)$.
If we don't wish to distinguish between the two mutually exclusive
possibilities of Definition~1 (ii), (iii), then wide-$(s)$ sequences are
more appropriate.

Throughout, let $X,Y$ be Banach spaces.

\proclaim{Proposition 2}
\roster
\item"(a)" A bounded subset $W$ of $X$ is non-relatively weakly compact
iff $W$ contains a wide-$(s)$ sequence.
\item"(b)" $T\in \L (X,Y)$ is non-weakly compact iff there is a sequence
$(x_n)$ in $X$ with $(x_n)$ and $(Tx_n)$ both wide-$(s)$ sequences.
\endroster
\endproclaim

\demo{Proof}

(a) First suppose $W$ is non-relatively weakly compact.
Choose $(x_n)$ in $W$ with no weakly convergent subsequence.
If $(x_n)$ has a weak-Cauchy subsequence $(x'_n)$, then $(x'_n)$ is of
course non-trivial, hence  $(x'_n)$ has an $(s)$-sequence by Proposition~2.2
of \cite{R3}.
If $(x_n)$ has no weak-Cauchy subsequence, $(x_n)$ has an $\ell^1$-subsequence
by the $\ell^1$-theorem.

Conversely, no wide-$(s)$ sequence can have a weakly convergent subsequence,
by Proposition~1, thus proving (a).
Alternatively, rather than using Proposition~1, we may give the following
elementary argument.

It is trivial that any subsequence of a wide-$(s)$ sequence is also 
wide-$(s)$.

Thus suppose to the contrary, that $(x_n)$ is a wide-$(s)$ sequence with
$(x_n)$ converging weakly to some $x$.
We may assume without loss of generality that $(x_n)$ is a basis for $X$,
since $x$ is in the closed linear span of the $x_n$'s.
{\it Define $\un{s}$}, the summing functional in $X^*$ by
$$\un{s} \Big(\sum c_j x_j\Big) = \sum c_j\ .$$
Then since $\un{s} (x_j) = 1$ for all $j$, $\un{s} (x)=1$, so $x\ne0$.
However no basic sequence can converge weakly to a non-zero element.
Indeed, if $(x_j^*)$ denotes the sequence of functionals biorthogonal to
$(x_j)$, then
$$\align
&x_i^* (x) = \lim_{j\to\infty} x_i^* (x_j) =0 \ \text{ for all } i\ ;
\text{ since}\\
&x= \sum_{i=1}^\infty x_i^* (x)x_i \text{ because }(x_j)
\text{ is a {\it basis\/} for } x\ ,\ x=0\ .\tag $*$
\endalign$$
\enddemo

\demo{Proof of (b)}
First suppose $T$ is non-weakly compact.
Then  $W\eqdf T( BaX)$ is non-relatively weakly compact.
Choose then $(x_j)$ in $Ba\,X$ with $(Tx_j)$ a wide-$(s)$ sequence.
But then $(x_j)$ cannot have a weakly convergent subsequence $(x'_j)$,
for else  $(Tx'_j)$ would be wide-$(s)$ and weakly convergent which is
impossible.
Hence $(x_j)$ has a wide-$(s)$ subsequence $(x'_j)$ by part~(a); then
also $(Tx'_j)$ is wide-$(s)$.

Of course conversely if $(x_j)$ and $(Tx_j)$ are wide-$(s)$ then $(x_j)$
is a bounded sequence and $(Tx_j)$ has no weakly convergent subsequence,
so $T$ is not weakly compact.\qed 
\enddemo

\proclaim{Corollary}
A Banach space is non-reflexive iff it contains a wide-$(s)$ sequence.
\endproclaim

\remark{Remark}
In Theorem~12 below, we prove a general result for selecting wide-$(s)$
sequences, which immediately yields Proposition~2(a), and hence this
corollary, without using the $\ell^1$-Theorem.
\endremark

We next review some permanence properties of wide-$(s)$ sequences.
First, a companion definition.

\definition{Definition 2}
A seminormalized sequence $(e_j)$ in a Banach space $X$ is called
\roster
\item"(i)" a wide-$(c)$ sequence if $(e_j)$ is a basic sequence  with
bounded partial sums; i.e., $\sup_n \|\sum_{j=1}^n e_j\| <\infty$.
\item"(ii)" a $(c)$-sequence if $(e_j)$ is a basic sequence with
$(\sum_{j=1}^n e_j)_{n=1}^\infty$ a weak-Cauchy sequence.
\item"(iii)" the difference sequence of a sequence $(b_j)$ if
$e_j = b_j-b_{j-1}$ for all $j>1$, $e_1 = b_1$.
\endroster
\enddefinition

\remark{Remark}
In the definition of ``wide-$(s)$'', the term ``$(s)$'' stands for
``summing.''
In the above, ``$(c)$'' stands for ``convergent''; of course here, the
series $\sum e_j$ is only weak Cauchy convergent; i.e., $\sum_j x^* (e_j)$
converges for all $x^*\in X^*$.
\endremark

As noted in the introduction, wide-$(s)$ and wide-$(c)$ sequences were
originally introduced by I.~Singer \cite{S1}.
He used the terminology ``$(x_j)$ is of type $P$ (resp. type $P^*$)''
if $(x_j)$ is a wide-$(c)$ (resp. wide-$(s)$) basis for a Banach space $X$.
The following result summarizes various permanence properties.
(For the proof of Proposition 3(i)--(iii), see \cite{S1} or Theorem~9.2
of \cite{S2}; Proposition~3(iv) is given as Proposition~2.1 of \cite{R3}.)

\proclaim{Proposition 3}
Let $(b_j)$ be a given sequence in $X$ with difference sequence $(e_j)$.
\roster
\item"(i)" $(b_j)$ is wide-$(s)$ iff $(e_j)$ is wide-$(c)$.
\item"(ii)" $(b_j)$ is wide-$(s)$ iff $(b_j)$ is bounded, $(\|e_j\|)$ is
bounded below, and $(e_j)$ is a basic sequence.
\item"(iii)" Assuming $(b_j)$ is a basic sequence with biorthogonal
functionals $(b_j^*)$ in $[b_j]^*$, then $(b_j)$ is wide-$(s)$ iff
$(b_j^*)$ is wide-$(c)$ iff $(\sum_{j=1}^n b_j^*)_{n=1}^\infty$ is wide-$(s)$.
\item"(iv)" $(b_j)$ is $(s)$ iff $(e_j)$ is $(c)$.
\item"(v)" If $(b_j)$ is a convex block-basis of a wide-$(s)$ sequence,
then $(b_j)$ is wide-$(s)$.
\endroster
\endproclaim

\remark{Remarks}
1. Proposition 3 (ii) yields that wide-$(s)$ sequences are 
{\it characterized\/} as semi-normalized basic sequences whose difference
sequence is also a basic sequence.

2. The statements (ii)--(iv) of Proposition 3 may also be deduced from
some arguments in \cite{R3} (specifically, see the proofs of Propositions~2.1
and 2.4).
\endremark

For the sake of completeness, we give the proof of Proposition~3(v).
Suppose then $(u_i)$ is a wide-$(s)$ sequence and there exist $0\le n_1 <
n_2 <\cdots$ and numbers $\lambda_1,\lambda_2,\ldots$ with
$$b_j = \sum_{i=n_j+1}^{n_{j+1}} \lambda_i u_i\ ,\quad
\lambda_i \ge0\ ,\quad
\sum_{i=n_j+1}^{n_{j+1}} \lambda_i =1\ \text{ for all } j\ .$$
Of course then $(b_j)$ is a basic sequence, since it is a block basis of one.
Since $(u_i)$ is wide-$(s)$, there is a number $\beta$ so that
$$\Big| \sum_{i=1}^k c_i\Big| \le \beta \| \sum c_i u_i\| \text{ whenever }
\sum c_i u_i\text{ converges.}$$
Now suppose scalars $c_1,c_2,\ldots$ given with only finitely many non-zero
and let $k\ge1$.
But then
$$\eqalignno{
\Big| \sum_{j=1}^k c_j\Big|
& = \Big| \sum_{j=1}^k c_j \sum_{i=n_j+1}^{n_{j+1}} \lambda_i\Big|\cr
& = \Big| \sum_{j=1}^k \sum_{i=n_{j+1}}^{n_{j+1}} c_j\lambda_i\Big|\cr
&\le \beta \Big\| \sum_j c_j \sum_{i=n_j+1}^{n_{j+1}} \lambda_i u_i\Big\|\cr
& = \beta\Big\| \sum c_j b_j\Big\| \ .&\qed\cr}$$

Of course Propositions~2 and 3 yield the following immediate consequence.

\proclaim{Corollary 4}
$T\in\L(X,Y)$ is non-weakly compact iff there is a sequence $(e_j)$ in $X$
so that $(e_j)$ and $(Te_j)$ are both wide-$(c)$ sequences.
\endproclaim

For the next result, recall that $T\in \L(X,Y)$ is {\it Tauberian\/} if
$T^{**}(X^{**}\sim X)\subset Y^{**} \sim Y$.

\proclaim{Theorem 5}
Let $T\in \L(X,Y)$.
The following are equivalent.
\roster
\item"(i)" $T$ is non-Tauberian.
\item"(ii)" There exists a wide-$(s)$ sequence $(x_j)$ in $X$ with
$(Tx_j)$ norm-convergent.
\item"(iii)" There exists a wide-$(c)$ sequence $(e_j)$ in $X$ with
$\sum \|Te_j\| <\infty$.
\endroster
\endproclaim

\proclaim{Corollary 6}
Again let $T\in \L(X,Y)$.
The following are equivalent.
\roster
\item"(i)" $T$ is Tauberian.
\item"(ii)" For every wide-$(s)$ sequence $(x_j)$ in $X$, there is a
subsequence $(x'_j)$ with $(Tx'_j)$ wide-$(s)$.
\item"(iii)" For every wide-$(c)$ sequence $(e_j)$, there exist $n_1<n_2<
\cdots$ with $(\sum_{j=n_i+1}^{n_{i+1}} Te_j)_{i=1}^\infty$ a wide-$(c)$
sequence.
\endroster
\endproclaim

\demo{Proof of Theorem 5}
Assume first $T$ is non-Tauberian.
it follows that there exists an $x^{**} \in X^{**} \sim X$ with $\|x^{**}\|=1$
and $T^{**}x^{**} \eqdf y$  in $Y$.
We shall prove there exists a normalized wide-$(s)$ sequence $(x_j)$ in $X$
with $(Tx_j)$ converging in norm to $y$.
It is convenient to isolate the following step of the proof.
\enddemo

\proclaim{Lemma 7}
Given $\ep >0$ and $F$ a finite dimensional subspace of $X^*$, there exists
an $x\in X$ with $\|x\| \le 1+\ep$, $\|Tx-y\|<\ep$, and
$$f(x) = x^{**} (f)\text{ for all } f\in F\ .
\tag 1$$
\endproclaim

\demo{Proof}
It is a standard result that given $G$ a finite-dimensional subspace of
$X^*$, there exists and $x_G$ in $X$ with $\|x_G\|< 1+\ep$ and
$$g(x) = x^{**} (g)\text{ for all } G\ .
\tag 2$$
But then letting $\D$ be the directed set of all finite-dimensional
subspaces of $X^*$ containing $F$, we obtain a net $(x_G)_{G\in \D}$ with
$\lim_{G\in \D} x_G = x^{**}$ weak*, with $x_G$ satisfying (1) for all
such $G$.
(Here, we regard $X\subset X^{**}$.)
Hence $\lim_{G\in\D} Tx_G =y$ weakly.
But then there exists  $x$ a convex combination of the $x_G$'s with
$\|Tx-y\| <\ep$.
But then $x$ satisfies (1) since every $x_G$ does.\qed
\enddemo

Now continuing with the proof of Theorem 5, since $x^{**}\notin X$,
$x^{**\bot} \eqdf \{x^*\in X^* :x^{**} (x^*)=0\}$ $\delta$-norms $X$
for some $\delta\ge 0$; i.e.,
$$\text{for all }x\in X\ ,\ \delta\|x\| \le \sup \{|z(x)|: z\in x^{**\bot} \ ,
\ \|z\|\le 1\}
\tag 3$$
Next, choose $x^*\in X^*$ with $x^{**}(x^*)=1$.

Now applying Lemma 8, we may inductively choose a sequence $(x_n)$ in $X$
and a sequence $(F_n)$ of finite-dimensional subspaces of $x^{**\bot}$
satisfying the following properties for all $n$ (where $F_0 =\{0\}$).
$$\align
& 1- {1\over2^n}  \le \|x_n\| \le 1+ {1\over 2^n} \tag 4\\
& \|Tx_n-y\| < {1\over 2^n}\tag 5\\
& F_n \ 2\delta\text{-norms}[x_1,\ldots,x_n] \tag 6\\
& F_{n-1} \perp x_n \tag 7\\
& x^* (x_n) = 1\tag 8\\
& F_{n-1} \subset F_n\ .\tag 9
\endalign$$
Indeed, suppose $n\ge1$ and $F_{n-1}$ chosen.
Applying Lemma~8 (and its obvious consequence that for an appropriate $F$,
any $x$ satisfying its conclusion must satisfy $\|x\|\ge 1-\ep$), we
let $F=F_{n-1} +[x^*]$; then choose $x=x_n$ satisfying (1) with $\|x\| \ge
1-{1\over2^n}$.
Since (1) holds, $x^{**}(f) = f(x) =0$ for all $f\in F_{n-1}$; i.e., (7) holds.

Of course also (8) holds.
Now since $x^{**\bot}$ $\delta$-norms $X$, we may choose $F_n\supset F_{n-1}$
a finite dimensional subspace of $x^{**\bot}$ which $2\delta$-norms
$[x_1,\ldots,x_n]$.

This completes the inductive construction.
Now if we let $Z= \bigcup_{i=1}^\infty F_n$, then (6) yields that
$$Z\ \ 2\delta\text{-norms} [x_j]_{j=1}^\infty\ .
\tag 10$$
Evidently (5) yields then $Tx_n\to y$.
Now we claim
$$(x_j) \text{ has no weakly convergent subsequence.}
\tag 11$$

Indeed, if not, let $(x'_j)$ be a subsequence weakly convergent to $x$, say.
Then by (8), $x^*(x) =1$.
But by (7), $f(x) = \lim_{n\to\infty} f(x_n)=0$ for all $f\in Z$.
Since $Z$ $2\delta$-norms $[x_j]_{j=1}^\infty$, $x=0$, a contradiction.

Since (11) is proved, $(x_j)$ has a wide-$(s)$ subsequence $(x'_j)$ by
Proposition~2(a).
Of course then $(x'_j/\|x'_j\|)$ is the desired normalized wide-$(s)$ sequence
whose image tends to $y$.
Thus (i) $\Rightarrow$ (ii) is proved.

(ii) $\Rightarrow$ (iii).
By passing to a subsequence, assume $\|Tx_j-Tx_{j-1}\|< {1\over2^j}$.
Then the difference sequence $(e_j)$ of $(x_j)$ satisfies (iii).

(iii) $\Rightarrow$ (ii).
Let $b_n = \sum_{j=1}^n e_j$ for all $n$.
Then $\lim_{j\to\infty} Tb_j = \sum_{j=1}^\infty Te_j$ in norm and $(b_j)$ is
wide-$(s)$.

(ii) $\Rightarrow$ (i).
Since $(x_j)$ is wide-$(s)$, it follows (e.g., by applying Proposition~2(a)),
that there is a weak*-cluster point $x^{**}$ of $(x_j)$ in $X^{**}\sim X$.
But then if $Tx_j\to y$ in $X$ say, $T^{**}x^{**}=y$, so $T$ is
non-Tauberian.\qed

\remark{Remarks}
1. Note the following  immediate consequence of Theorem~5:
$T$ is Tauberian iff $T|Z$ is Tauberian for every separable $Z\subset X$
iff $T|Z$ is Tauberian for every $Z\subset X$ with a basis.

2. We prove a stronger quantitative version (called Theorem~$5'$) later on,
in the last remark following the proof of Corollary~15.

3. As proved in \cite{NR}, if $T$ is non-Tauberian, there exists $Z$ a
closed linear subspace of $X$ with $T(Ba\,Z)$ non-closed.
If we then choose $y\in \overline{T(Ba\,Z)} \sim T(Ba\,Z)$, we may choose
$(x_n)$ in $Ba(Z)$ with $Tx_n\to y$. 
It now follows that any weak*-cluster point $G$ of $(x_n)$ in $Z^{**}$
does not belong to $Z$.
Of course then $(x_n)$ has no weakly convergent subsequence; it follows
then by known results that $(x_n)$ has either an
$(s)$-subsequence or an $\ell^1$-subsequence.
We preferred here, however, to give a direct self-contained proof of the
non-trivial implication in Theorem~5.
\endremark

\demo{Proof of Corollary 6}
(i) $\Rightarrow$ (ii).
Let  $(x_j)$ be a wide-$(s)$ sequence in $X$.
Then it follows that $(Tx_j)$ has no weakly convergent subsequence.
Indeed, otherwise,  there would exist $(b_j)$ a convex block basis of $(x_j)$
with $(Tb_j)$ norm-convergent.
But then $(b_j)$ is wide-$(s)$ by Proposition~3(v), hence $T$ is
non-Tauberian by Theorem~5, a contradiction.
Thus by Proposition~2(a), there exists $(x'_j)$ a subsequence of $(x_j)$
(with $(Tx_j)$ wide-$(s)$; of course $(x'_j)$ is still wide-$(s)$, so this
implication is proved.

(ii) $\Leftrightarrow$ (iii)
is immediate from the permanence property Proposition~3(i) and the evident
fact that if $(b_j)$ is a given sequence and $0\le n_0<n_1<\cdots$, then if
$(e_j)$ is the difference sequence of $(b_j)$,
$(\sum_{j=n_i+1}^{n_{i+1}} e_j)_{i=1}^\infty$ is the difference sequence
of $(b_{n_i})_{i=1}^\infty$.
Of course no wide-$(s)$ sequence can be norm-convergent, so (ii) $\Rightarrow$
(i) follows immediately from Theorem~5.\qed
\enddemo

Theorem 5 easily yields the following result:

\proclaim{Corollary 8}
Let $T\in \L(X,Y)$. The following are equivalent\/:
\roster
\item"(a)" {\it $T$ is Tauberian}
\item"(b)" {\it $\ker (T+K)$ is reflexive for all compact $K\in \L(X,Y)$.}
\item"(c)" {\it $\ker (T+K)$ is reflexive for all nuclear $K\in \L(X,Y)$.}
\endroster
{\rm (The equivalence (a) $\Leftrightarrow$ (b) is due to Gonz\'alez and Onieva
\cite{GO}.) 
Indeed, (a) $\Rightarrow$ (b) is immediate since then $T+K$ is also
Tauberian, and (b) $\Rightarrow$ (c) is trivial.
Now suppose (c) holds yet (a) is false.
Then by Theorem~5(ii), after passing to a subsequence, we can choose
$(x_j)$ a wide-$(s)$ sequence in $X$ with $\|Tx_{j+1} - Tx_j\|< {1\over2^j}$
for all $j$.
It then  follows immediately that $T|[x_j]$ is nuclear.
Indeed, letting $(b_j)$ be the difference sequence for $(x_j)$, then $(b_j)$
is also a semi-normalized basis for $[x_j]$ and $\sum\|Tb_j\|< \infty$.
But then there exists a nuclear operator $K:X\to Y$ with $K|[x_j] = T|[x_j]$.
Hence $\ker (T-K)\supset [x_j]$ is non-reflexive, because $(x_j)$ is a
wide-$(s)$-sequence lying inside this kernel.}\qed 
\endproclaim

We next show that wide-$(s)$ sequences immediately yield the triangular
arrays in non-reflexive Banach spaces, discovered by R.C.~James \cite{J1},
\cite{J2}. 

\proclaim{Proposition 9}
Let $(b_j)$ be a wide-$(s)$ sequence in a Banach space $X$.
There exist bounded sequences $(f_j)$ and $(g_j)$ in $X^*$ satisfying
$$\align
\text{\rm (i)}&\qquad f_i(b_j)=1\text{ for all } j\ge i\\
&\qquad f_i(b_j)=0\text{ for all } j<i\\
\noalign{\vskip6pt}
\text{\rm (ii)}&\qquad g_j(b_i)=1\text{ for all } j\ge i\\
&\qquad g_j(b_i)=0\text{ for all } j<i\ .
\endalign$$
\endproclaim

\demo{Proof}
We may assume without loss of generality, by the Hahn-Banach Theorem, that
$X= [b_j]$. Let $\bs$ be the summing functional for $(b_j)$ and $(b_j^*)$
be the sequence biorthogonal to $(b_j)$.
Then simply let $g_j = \sum_{i\le j} b_i^*$ and $f_j =\bs-\sum_{i<j} b_i^*$
for all $j$.
Since $(b_j)$ dominates the summing basis, $(g_j)$ and hence $(f_j)$
are uniformly bounded.\qed
\enddemo

\remark{Remark}
It follows immediately from the duality theory given in the Proposition,
page 722 of \cite{R3}, that if $X= [b_j]$ as above, both sequences $(f_j)$ and
$(g_j)$ are again basic; in fact both are wide-$(s)$ sequences in $X^*$.
\endremark

We now give some quantitative definitions, in order to obtain further
permanence properties and obtain certain localizations of the preceding
results.

(Recall: If $(b_j)$ is a basic-sequence or a finite sequence, then $(b_j)$
is called a $\lambda$-basic sequence if $\|\sum_{j=1}^k c_j b_j\|
\le \lambda\|\sum c_j b_j\|$ for all $k$ and scalars $c_1,c_2,\ldots$ with
$\sum c_j b_j$ convergent.)

\definition{Definition 3}
A (finite or infinite) sequence $(b_j)$ in Banach space is called
$\lambda$-wide-$(s)$ if
\roster
\item"(a)" $(b_j)$ is a $2\lambda$-basic sequence. 
\item"(b)" $\|b_j||\le \lambda$ for all $j$.
\item"(c)" $|\sum_{j=k}^n c_j| \le \lambda \| \sum_{j=1}^n  c_j b_j\|$
for all $1\le k\le n<\infty $ (resp. with $n$ the length of $(b_j)$
if finite), and scalars $c_1,c_2,\ldots,c_n$. 
\endroster
\enddefinition 

\definition{Definition 4}
A (finite or infinite) sequence $(e_j)$ in a Banach space is called
$\lambda$-wide-$(c)$ if
\roster
\item"(a)" $(e_j)$ is a $\lambda$-basic sequence.
\item"(b)" $\|e_j\| \ge {1\over\lambda}$ for all $j$.
\item"(c)" $\|\sum_{j=1}^k e_j\|\le \lambda$ for all $k$.
\endroster
\enddefinition

\remark{Remarks}
1. Of course an infinite sequence $(b_j)$ (resp. $(e_j)$) is wide-$(s)$
(resp. wide-$(c)$) iff it is $\lambda$-wide-$(s)$ (resp. $\lambda$-wide-$(c)$)
for some $\lambda\ge1$.
Also, note that if $(b_j)$ is $\lambda$-wide-$(s)$, then
trivially $\|b_j\| \ge {1\over\lambda}$ and
$\|b_j^*\| \le 2\lambda$  for all $j$, where $(b_j^*)$ is biorthogonal
to $(b_j)$ in $(b_j)^*$.
Also immediately $(e_j)$ $\lambda$-wide-$(c)$ implies $\|e_j\|\le2\lambda$
for all $j$.

2. An easy refinement of the (easy) proof of Proposition~3(v) shows that any
convex block-basis of a (finite or infinite) $\lambda$-wide-$(s)$ sequence
is again $\lambda$-wide-$(s)$.

3. It is easily seen that if $(b_i)$ is a wide-$(s)$ sequence in $c_0$,
then the basis-constant for $(b_i)$  is at least~2.
This is why we give the requirement that $\lambda$-wide-$(s)$ sequences
be $2\lambda$-basic.
\endremark

The following quantitative version of Proposition~3(i) follows immediately
from the arguments in \cite{R3}.

\proclaim{Proposition 10}
Let $(b_j)$ be a finite or infinite sequence with difference sequence $(e_j)$.
For all $\lambda\ge1$, there exists a $\beta\ge1$ so that
\roster
\item"(a)" If $(b_j)$ is $\lambda$-wide-$(s)$, then $(e_j)$ is
$\beta$-wide-$(c)$.
\item"(b)" If $(e_j)$ is $\lambda$-wide-$(c)$, then $(b_j)$ is
$\beta$-wide-$(s)$.
\endroster
\endproclaim

The next result, immediately yields the quantitative triangular array
result of R.C.~James \cite{J1}, in virtue of Proposition~9.

\proclaim{Theorem 11}
Every non-reflexive Banach space contains for every $\ep>0$, a
normalized $(1+\ep)$-wide-$(s)$ sequence.
\endproclaim

To deduce the quantitative result in \cite{J1}, suppose $(b_i)$ is
normalized and $(1+\ep)$-wide-$(s)$.
Let $f_j = \sum_{i=j}^\infty b_i^* = \bs- \sum_{i<j} b_i^*$ as in
Proposition~9.
Thus $\|f_j\| \le 1+\ep$ for all $j$.
Now set $h_j = f_j/\|f_j\|$ for all $j$.
Thus $(b_j)$ and $(h_j)$ are norm-one sequences satisfying
$$\eqalign{
&h_i (b_j) \ge {1\over 1+\ep} \ \text{ for all }\ j\ge i\cr
&h_i (b_j) = 0\text{ for all }\ j<i\ .\cr}
\tag 12$$

\remark{Remark}
The work in \cite{J2} (specifically Theorem 8) essentially yields
Theorem~11.
(I had overlooked the fundamental reference \cite{J2} in the earlier
version of this paper.)
The treatment given here provides a general criterion for selecting
$(1+\ep)$-wide-$(s)$ sequences out of a given sequence, which yields the
apparently stronger result Theorem~11$'$ below.
\endremark

Theorem 11 is proved by a refinement of the argument given for Proposition~2.2
in \cite{R3}.
We will in fact show the following result, which easily yields Theorem~11.
(As usual, we regard $X\subset X^{**}$; for $G$ in $X^{**}$,
$\dist (G,X)$ denotes the distance of $G$ to $X$; i.e., $\dist (G,X) =
\inf \{\|G-x\|: x\in X\}$.

\proclaim{Theorem 12}
Let $(x_j)$ be a bounded sequence in a Banach space $X$, having a
weak*-cluster point $G$ in $X^{**}\sim X$.
Let $d=\dist (G,X)$, $\ep>0$, and $\lambda = {1\over d} +\ep$,
$\beta = {\|G\|+d\over d} +\ep$.
Then $(x_j)$ has a subsequence $(b_j)$ satisfying the following:
\roster
\item"1)" $(b_j)$ is $\beta$-basic.
\item"2)" $|\sum_{j=k}^n c_j| \le \lambda \|\sum_{j=1}^n c_jb_j\|$
for all $1\le k\le n<\infty$ and scalars $c_1,\ldots,c_n$.
\endroster
\endproclaim

\proclaim{Corollary}
If
$$\gamma = \max \{\frac1d ,{\|G\|+d\over 2d},\sup_j \|b_j\|\}$$
then for any $\ep>0$, $(x_j)$ has a $(\gamma+\ep)$-wide-$(s)$-subsequence.
\endproclaim

Of course Theorem 11 follows immediately from this result.
Indeed, we may assume without loss of generality that $X$ is separable
non-reflexive.

Let $\ep>0$, and let $\delta>0$ be decided.
Choose (using Riesz's famous lemma) a $G$ in $X^{**}$ with $\|G\|=1$
and $\dist (G,X) >1-\delta$.
Next choose $(x_j)$ a normalized sequence in $X$ having $G$ as a
$\omega^*$-cluster point.
Now if $\delta$ is such that $\frac1{1-\delta} <1+\ep$, then $(x_j)$
has a $1+\ep$-wide-$(s)$ subsequence by the Corollary.

For the proof of Theorem 12 we first recall some standard ideas and
results (cf.\ \cite{R3}).
Given $0<\eta \le1$ and $Y$ a linear subspace of $X^*$, $Y$ is said to
$\eta$-norm $X$ if
$$\eta \|x\| \le \sup_{y\in Ba\, Y} |y(x)|\ \text{ for all }\ x\in X\ .  $$

The next two lemmas summarize well known material.

\proclaim{Lemma 13}
Let $(x_j)$ be a semi-normalized sequence in $X$ and $Y$ an $\eta$-norming
subspace of $X^*$.
Assume that $y(x_j)\to0$ as $j\to\infty$ for all $y\in Y$.
Then given $0<\ep<\eta$, $(x_j)$ has a $\frac1{\eta-\ep}$-basic subsequence.
\endproclaim

\proclaim{Lemma 14}
Let $G\in X^{**}\sim X$, $\delta =\dist (G/\|G\|,X)$,
$G_\bot = \{x^*\in X^* : G(x^*) =0\}$.
Then $G_\bot$ $\frac{1+\delta}\delta$-norms $X$.
\endproclaim

\remark{Remark}
Lemmas 13 and 14 immediately imply the classical result of M.I.~Kadec and
A.~Pe{\l}czy\'nski \cite{KP} that a semi-normalized sequence in a Banach
space has a basic subsequence provided every weakly convergent subsequence
converges weakly to zero.
These lemmas are in fact a crystallization of the arguments in \cite{KP}.
\endremark

\demo{Proof of Theorem 12}
We have that $\dist (\frac{G}{\|G\|} ,X) = \frac{d}{\|G\|}$.
It follows immediately from Lemmas~13 and 14 that $(x_j)$ has a
$\beta$-basic subsequence, so let us assume then that $(x_j)$ is
already $\beta$-basic.
Let $\eta >0$ be decided later.
We shall choose $(b_j)$ a subsequence of $(x_j)$ and a sequence $(f_j)$
in $X^*$ satisfying the following conditions for all $n$:
$$\align
&\|f_n\| < \frac1d +\eta \tag 13\\
&f_i(b_j) = 0\ \text{ for all }\ 1\le j<i\le n \tag 14\\
&G(f_i)=1\ \text{ for all }\ 1\le i\le n\tag 15\\
&|f_i(b_n)-1| < {\eta\over 2^n}\ \text{ for all }\ 1\le i\le n\ .\tag 16
\endalign$$

First, by the Hahn-Banach theorem, choose $F\in X^{***}$ with
$$\|F\| = \frac1d\ ,\quad F(G)=1\ ,\ \text{ and }\
F(x)=0\ \text{ all }\ x\in X\ .
\tag 17$$
Now choose $f_1$ in $X^*$ satisfying (13) and (15) for $n=1$.
Suppose $n\ge1$ and $f_1,\ldots,f_n$, $b_1,\ldots, b_{n-1}$ chosen
satisfying (13)--(15); suppose also $b_{n-1} = x_k$ (if $n>1$).
Since (15) holds and $G$ is a $w^*$-cluster point of the $x_j$'s, we may
choose an $\ell>k$ so that setting $b_n= x_\ell$, then (16) holds.
Since the span of $b_1,\ldots, b_n$ and $G$, 
$[b_1,\ldots,b_n,G]$, is
finite dimensional, by (17) we may choose $f_{n+1}$ in $X^*$ with
$\|f_{n+1}\| <\frac1d +\eta$, so that $f_{n+1}$ agrees with $F$ on
$[b_1,\ldots,b_n,G]$; i.e., $G(f_{n+1})=1$ and  $f_{n+1}(b_i)=0$
all $1\le i\le n$.
This completes the inductive construction of $(b_j)$ and $(f_j)$.

Now let $1\le k\le n$ and scalars $c_1,\ldots,c_n$ be given with
$\|\sum_{j=1}^n c_jb_j\|\le 1$.
$$\align
\Big| \sum_{j=k}^n c_j\Big|
& = \Big| f_k\biggl( \sum_{j=1}^n c_jb_j\biggr) +
\sum_{j=k}^n c_j (1-f_k(b_j))\Big|\ \text{ by (14)}\\
&\le {1\over d} +\eta + \sup_j |c_j|\eta\ \text{ by (13) and (16)}\\
&\le {1\over d} +\eta +\sup_j \|b_j^*\| \eta\
\text{ (where $(b_j^*)$ is the sequence biorthogonal to $(b_j)$ in
$[b_j]^*$)}\ .
\endalign$$
Since $(x_j)$ was assumed basic, so is $(b_j)$, whence $\tau \eqdf \sup_j
\|b_j^*\| <\infty$.
Thus if $\eta$ is such that $\eta (1+\tau)\le\ep$, the proof is finished.\qed
\enddemo

The next result yields absolute constants in the selection of general basic
sequences in Banach spaces.

\proclaim{Corollary 15}
Let $(x_j)$ be a bounded non relatively compact sequence
in the Banach space $X$.
Given $\ep >0$, there exists an $x\in X$ and a subsequence $(x'_j)$ of
$(x_j)$ so that the following holds, where $b_j=x'_j-x$ for all $j$.
\roster
\item"1)" If $\{x_1,x_2,\ldots\}$ is relatively weakly compact, then
$(b_j)$ is a $(1+\ep)$-basic sequence.
\item"2)" If $\{x_1,x_2,\ldots\}$ is not relatively weakly compact, then
$(b_j)$ is a $(2+\ep)$-basic sequence.
Moreover there exists a block basis $(u_j)$ of $(b_j)$ and a $c>0$ so
that $(cu_j)$ is a $(1+\ep)$-wide-$(s)$ sequence.
\endroster
\endproclaim

\demo{Proof}
Since $(x_j)$ is non-relatively compact, we may by passing to a subsequence
assume that for some $\delta >0$, $\|x_i-x_j\| \ge\delta$ for all $i\ne j$.

In Case~1), choose $(\bar x_j)$ a subsequence of $(x_j)$ and $x\in X$
with $\bar x_j\to x$ weakly.
Then $(\bar x_j-x)_{j=1}^\infty$ satisfies the hypotheses of Lemma~13 for
$\eta=1$, where ``$Y$''~$=X^*$. 
The conclusion of 1) is now immediate from Lemma~13.

In Case 2), we may choose $G\in X^{**} \sim X$ with $G$ a $w^*$-cluster point
of $(x_j)$.
Now let $d=\dist (G,X^{**})$, $\eta >0$ to be decided later, and
choose $x\in X$ with
$$\|G-x\| < d+\eta\ .
\tag 18$$
Now setting $H=G-x$, we have by Riesz's famous argument that since
$$\align
&d\le \|H\| \le d+\eta\ ,\tag 19\\
&\dist \left( {H\over \|H\|} ,X\right) \ge {d\over d+\eta}\ .\tag 20
\endalign$$
Thus if $\frac{\eta}d <\ep$, it follows from Lemmas~13, 14, and the fact
that $H$ is a weak*-cluster point of $(x_j-x)$, that $(x_j)$ has a
subsequence $(x'_j)$ satisfying the first statement in (2).
Finally, it follows by the techniques in \cite{R1} that there exists a
convex block basis $(u_j)$ of $(b_j)$ and a separable isometrically
norming subspace $Y$ of $X^*$ so that
$$\|u_j\| <d+\eta\ \text{ for all }\ j
\tag 21$$
and
$$y(u_j)\to H(y) \text{ as } j\to \infty\ ,\ \text{ all } y\in Y \ .
\tag 22$$
Now (19), (21) and (22) yield that if $F$ is any $w^*$-cluster point
of $(u_j)$ in $X^{**}$, then
$$\|F\| \le d+\eta\ \text{ and }\ \dist (F,X)\ge d\ .
\tag 23$$
Fixing such an $F$, let $c= \|F\|^{-1}$; then of course $(cu_j)$ has
$F/\|F\|$ as a $w^*$-cluster point.
It now follows from the Corollary to Theorem~12 (still assuming $\eta<\ep d$)
that $(cu_j)$ has a $(1+\ep)$-wide-$(s)$ subsequence.\qed
\enddemo

\demo{Remarks}
1. We may now easily obtain the following quantitative refinement of
the non-trivial part of Theorem~5.
\enddemo

\proclaim{Theorem 5$'$}
If $T\in \L(X,Y)$  is non-Tauberian, then for all $\ep>0$, there exists
a $(1+\ep)$-wide-$(s)$ sequence $(b_j)$ in $X$ with $(Tb_j)$ norm-convergent. 
\endproclaim

\demo{Proof}
We may deduce this by just quantifying our proof of Theorem~5 and
Corollary~15. 
We prefer, however, to deduce the result directly from Theorem~5 and
Corollary~15.
Choose $(x_j)$ wide-$(s)$ with $(Tx_j)$ norm-convergent,
by Theorem~5.
Now choose by Corollary~15, a convex block basis $(u_j)$ of $(x_j)$,
an element $x\in X$, and a $c>0$ so that $b_j \eqdf c(u_j -x)$ is
$1+\ep$-wide-$(s)$.
But of course since $(Tx_j)$ is norm-convergent, so is $(Tb_j)$.\qed
\enddemo

2.  We also obtain a quantitative refinement of Proposition~2(a), which
moreover sharpens Theorem~11, namely

\proclaim{Theorem 11$'$}
If $T\in \L(X,Y)$ is non-weakly compact, then for all $\ep>0$, there
exists  a $(1+\ep)$-wide-$(s)$ sequence $(b_n)$ in $X$ with $(Tb_n)$
wide-$(s)$.
\endproclaim

\demo{Proof}
By Proposition 2(a), first choose $(x_j)$ in $X$ with $(x_j)$ and $(Tx_j)$
wide-$(s)$.
By Corollary~15, $\ep>0$ given, there is a convex block basis $(u_j)$ of
($x_j)$, a $u$ in $X$, and a $c>0$, so that $(b_j)$ is $(1+\ep)$-wide-$(s)$,
where $b_j= c(u_j-u)$ for all $j$.
But of course $(Tu_j)$ is again wide-$(s)$, so $(Tb_j)$ cannot have a
weakly convergent subsequence, hence finally we can choose a
subsequence $(b'_j)$ of $(b_j)$ with $(Tb'_j)$ also wide-$(s)$.\qed

3. It is easily seen that the  $(1+\ep)$-wide-$(s)$ sequences in
Theorem~$5'$ and Theorem~$11'$  may be chosen to be normalized.
In fact, in most of the results formulated here, we can take our sequences
normalized.
Indeed, suppose $(b_j)$ is an arbitrary semi-normalized sequence in a Banach 
space.
Then we may choose $b'_j$ a subsequence  and a positive number $c$ so that
$\|b'_j\|$ rapidly converges to $c$, say $\big|\,\|b'_j\|-c\big|<\frac1{2^j}$ 
for all $j$.
Then $\eta >0$ given, it follows from standard perturbation results that
for some $N$, the sequence $(b'_j)_{j=N}^\infty$ is $1+\eta$-equivalent to
$(cb'_j/\|b'_j\|)_{j=N}^\infty$.
Now if $(b_j)$ is $\lambda$-wide-$(s)$, then of course ${1\over\lambda}\le
c\le \lambda$ and we obtain that $(b'_j/\|b'_j\|)_{j=N}^\infty$ is
$(c\vee {1\over c}) (1+\eta)^2\lambda$-wide-$(s)$.
So of course this normalized sequence is $\lambda^2(1+\eta)^2$-wide-$(s)$.
Thus if $\ep>0$ is given and $\lambda^2 (1+\eta)^2 <1+\ep$, $\lambda>1$,
and $(b_j)$ is $\lambda$-wide-$(s)$ with $(Tb_j)$ norm-convergent, then
also $(Tu_j)$ is norm-convergent, where $u_j=b'_{j+N}/\| b'_{j+N}\|$
all $j$, and $(u_j)$ is $(1+\ep)$-wide-$(s)$.

4. For a uniformity estimate in the biorthogonal constant of uniformly
separated bounded sequences, see \cite{HKPTZ}.
\enddemo

We next pass to the stability of wide-$(s)$ sequences under perturbations.
We first show that after passing to subsequences, triangular arrays and
wide-$(s)$ sequences are manifestations of the same phenomena.

\definition{Definition 5}
Given $\lambda\ge1$, a finite  or infinite sequence $(b_j)$ in a Banach
space is called $\lambda$-triangular if there exists a sequence $(f_j)$
in $X^*$ so that
$$\eqalign{
\text{\rm (i)}\qquad &f_i(b_j) = 1\ \text{ for all }\ j\ge i\cr
&f_i(b_j) = 0\ \text{ for all }\ j<i\cr
\noalign{\text{and}}
\text{\rm (ii)}\qquad &\|f_j\| ,\|b_j\|\le\lambda\ \text{ for all }\ j\ .\cr}$$
\enddefinition

Evidently Proposition 9 and Definition 3 yield immediately that
{\it every\/} (finite or infinite) {\it $\lambda$-wide-$(s)$ sequence
$(b_j)$ is $\lambda$-triangular\/}.
Thus Theorem~11 yields immediately that every non-reflexive Banach
space has, for every $\ep >0$, a normalized $(1+\ep)$-triangular sequence.
(This result is due to R.C.~James; see Theorem~8 of \cite{J2}.)
Now conversely, suppose $(b_j)$ is an infinite $\lambda$-triangular sequence,
and assume without loss of generality that $X= [b_j]$.
Now let $G$ be a $w^*$-cluster point of $(b_j)$ in $X^{**}$.
Then of course $\|G\| \le\lambda$, and $G(f_i)=1$ for all $i$.
Then letting $F$ be a $w^*$-cluster point of $(f_i)$ in $X^{***}$, it follows
that $F(G) = 1$ and $F(x) =0$ all $x\in X$.
Hence
$$\dist (G,X) \ge \frac1{\lambda}\ .
\tag 24$$
The Corollary to Theorem 12 now immediately yields

\proclaim{Corollary 16}
If $(b_j)$ is a $\lambda$-triangular sequence, then for every $\ep>0$,
$(b_j)$ has a $(\frac{\lambda^2+1}2 +\ep)$-wide-$(s)$ subsequence.
\endproclaim

A refinement of this reasoning now yields the following perturbation result.

\proclaim{Corollary 17}
Let $\lambda\ge1$ be given.
There exist $\beta\ge1$ and $\ep>0$ so that if $(b_j)$ and $(p_j)$ are 
infinite sequences in a Banach space with $(b_j)$ $\lambda$-wide-$(s)$ and
$\|p_j\|\le\ep$ for all $j$, then
$$\text{the sequence $(b_j+p_j)_{j=1}^\infty$ has a $\beta$-wide-$(s)$
subsequence.}
\tag 25$$
\endproclaim

\demo{Proof}
Let $\ep>0$; we shall discover the appropriate bounds for $\ep$ and $\beta$
in the course of the argument.
Choose $(f_j)$ in $X^*$ satisfying (i), (ii) of Definition~5.
Since $(b_j)$ is $\lambda$-triangular, it follows as in the argument preceding
Corollary~16 that if $G$ is any $w^*$-cluster point of $(b_j)$ in $X^{**}$,
then
$$\|G\|\le \lambda\ \text{ and }\ \dist (G,Y) \ge \frac1{\lambda}\ ,
\tag 26$$
where $Y=[b_j]$.
But then (cf.\ Lemma 2.6 of \cite{R2})
$$\dist (G,X) \ge \frac1{2\lambda}\ .
\tag 27$$
Of course if $P$ is any $w^*$-cluster point of $(p_j)$, then $\|P\|\le\ep$.
Hence we obtain that if $H$ is a weak*-cluster point of
$(b_j+p_j)_{j=1}^\infty$, then
$$\|H\| \le\lambda+\ep \ \text{ and }\ \dist (H,X)\ge \frac1{2\lambda}
-\ep\ .
\tag 28$$
Evidently we thus obtain from the Corollary to Theorem~12 that
$(b_j+p_j)_{j=1}^\infty$ has a $\beta$-wide-$(s)$ subsequence provided
$$\ep <\frac1{2\lambda}
\tag 29$$
where
$$\beta = \ep +  \Big( (\lambda^2 +\tfrac12) \vee 2\lambda\Big) \Big/
(1-2\ep\lambda)\ .
\tag 30$$
Evidently if we just let $\ep = \frac1{4\lambda}$, we obtain that
$(b_j+p_j)$ has a $(2\lambda^2 +5)$-wide-$(s)$ subsequence.\qed
\enddemo

\remark{Remarks}

1. If we assume the $p_j$'s lie in the closed linear span of the
$b_j$'s, we obtain the better estimate
$$\ep <\frac1{\lambda}\ \text{ and }\
\beta = {\lambda^2+1\over 2(1-\lambda\ep)} + \ep\ .$$
Indeed, this follows immediately from the Corollary to Theorem~12 and the
fact that then in (28), we have the better estimate $\dist (H,X)\ge
{1\over\lambda}-\ep$.
(Thus if $(b_j)$ is say $(1+\eta)$-wide-$(s)$,
$\eta$ small, then for $\ep$-sufficiently
small, the perturbed sequence  would have a $(1+2\eta)$-wide-$(s)$
subsequence.

2. Of course the proof only requires that $(b_j)$ is $\lambda$-triangular.
Moreover the same qualitative conclusion holds if we just assume
instead that $(p_j)$ is uniformly bounded so that all weak*-cluster points
have distance at most $\ep$ from $X$.
That is, we obtain then the following generalization.

3. Most of our results here are of a ``refined subsequence'' nature.
Nevertheless, wide-$(s)$ sequences and triangular arrays may have large spans.
For example, J.R.~Holub \cite{Ho} (cf.\ also \cite{S2}, pp.627--628) has
obtained that $L^1([0,1])$ has a wide-$(s)$ basis;
his argument yields also that
$C(\Delta)$ has a wide-$(s)$ basis, $\Delta$ the Cantor set, and hence $C(K)$
has such a basis, any uncountable compact metric space $K$.
\endremark

This suggests the following problem.

\proclaim{Question 1}
Does every non-reflexive Banach space with a basis have a wide-$(s)$
basis?
\endproclaim

Now triangular arrays may have large linear spans, even when the space has no
basis. This suggests

\proclaim{Question 2}
Does every separable non-reflexive Banach space have a $\lambda$-triangular
fundamental sequence for some $\lambda>1$? For every $\lambda >1$?
\endproclaim

\noindent (A sequence $(b_j)$ in a Banach space $B$ is called fundamental
if $[b_j]=B$.)

\proclaim{Corollary 17$'$}
Let $\lambda \ge1$, $\ep < {1\over2\lambda}$, and $(b_j)$, $(p_j)$ be
sequences in $X$ with $(b_j)$ $\lambda$-triangular, and $(p_j)$ bounded such
that all weak*-cluster points of $(p_j)$ in $X^{**}$ have distance at most
$\ep$ from $X$.
Then assuming $\|p_j\|\le M$ all $j$, there exists a $\beta$ depending only
on $\lambda,M$ and $\ep$, so that $(b_j+p_j)$ has a $\beta$-wide-$(s)$
subsequence.
\endproclaim

For the sake of definiteness, we note that the proof actually yields that
given $\eta>0$ (arbitrarily small), we may choose
$$\beta = {(\lambda^2 +\lambda M +\frac12 -\lambda\ep)\vee 2\lambda
\over 1-2\ep \lambda} + \eta\ .$$
Of course an interesting special case occurs when $(p_j)$ is weakly
convergent or even a constant sequence; then we have (since $\ep=0$) that
$$\beta = \left(\Big(\lambda^2 + \lambda M +\frac12\Big)
\vee 2\lambda\right) +\eta  \ ,$$
or in the case where the $p_j$'s lie in $[b_j]$,
$$\beta = \left(  (\lambda+M) \vee
\Big( {\lambda^2 +\lambda M+1\over 2}\Big) \right)
+\eta\ .$$

3. Combining the last observation in the preceding remark with Theorem~5$'$,
we obtain that {\it if $T\in \L(X,Y)$ is non-Tauberian, then for all
$\ep>0$ there exists $(b_j)$ a $(2+\ep)$-wide-$(s)$ sequence in $X$ with
$(Tb_j)$ norm-convergent and $\|Tb_j\| \le \ep$ for all $j$.}

Indeed, simply choose $(x_j)$ a $(1+\ep)$-wide-$(s)$ sequence in $X$
with $(Tx_j)$ norm-convergent, by Theorem~$5'$.
Now given $\ep>0$, choose $k$ with $\|Tx_n-Tx_k\| \le {\ep\over2}$ for
all $n\ge k$.
But then $(x_n-x_k)_{n=k+1}^\infty$ has a $\lambda$-wide-$(s)$
subsequence, with
$$\lambda \le \left( \Bigl( 1+{\ep\over2}\Bigr)^2 + \frac12 +\ep\right) \vee
(2+\ep) = 2+\ep \text{ for $\ep$ small enough.}
\eqno\qed$$
(We do not know if ``2'' can be replaced by ``1'' in the above assertion.)

We now localize some of our preceding results.

\proclaim{Proposition 18}
Given $\lambda\ge1$, $\ep>0$, and $k$, there is an $n$ so that every
$\lambda$-triangular sequence  of length at least $n$ has a
$({\lambda^2+1\over2}+\ep)$-wide-$(s)$ subsequence of length $k$.
\endproclaim

\demo{Proof}
We given an ``old fashioned'' compactness argument, just using
Corollary~16.
Were this false, we can find for every $n$ a norm $\|\cdot\|_n$ on
$\real^n$ so that if $(b_j)_{j=1}^n$ denotes the usual unit vectors,
then setting $X_n = (\real^n,\|\cdot\|_n)$, and letting
$f_j = \sum_{i=j}^n b_j^*$ for all $n$, we have that
$$\|b_j\|_n \le\lambda\ ,\ \|f_j\|_n^* \le\lambda\text{ all } 1\le j\le n\ ,
\tag 31$$
and $(b_j)_{j=1}^n$ has no subsequence of length $k$ which is
$\beta$-wide-$(s)$ in $X_n$, where $\beta = {\lambda^2+1\over2}+\ep$.

Now of course $(b_j^*)_{j=1}^n$, the functionals biorthogonal to $(b_j)$, 
are bounded by $2\lambda$.
But then it follows that regarding $(b_j)_{j=1}^\infty$ instead as the
usual unit basis of $c_{00}$, the set of all sequences which are ultimately
zero, we may choose $k_1<k_2<\cdots$ so that for all $x\in c_{00}$,
$$\lim_{n\to\infty} \|x\|_{k_n} \eqdf \|x\|\ \text{ exists.}
\tag 32$$
Of course this limit will exist uniformly on $W_n \eqdf \{x\in\real^n :
\|x\|_\infty \le1\}$, for each $n$.
Now let $X$ be the completion of $(c_{00},\|\cdot\|)$.
But then it follows immediately that $(b_j)_{j=1}^\infty$ is
$\lambda$-triangular in $X$, yet $(b_j)$ has no
$(\beta-{\ep\over2})$-wide-$(s)$ subsequence in $X$, of length $k$.
But by Corollary~16, $(b_j)$ {\it has\/} an infinite
$(\beta-{\ep\over2})$-wide-$(s)$ subsequence. 
This contradiction completes the proof.\qed
\enddemo

We next consider localized results which follow directly from our work above,
known facts about ultraproducts, and certain results in \cite{GM} to
which we refer for all unexplained concepts.
First, we briefly recall some ideas concerning ultraproducts.
Let $\U$ be a non-trivial ultrafilter on $N$.
$X_{\U}$, an ultraproduct of $X$, denotes the Banach space $\ell^\infty (X)/
\N_{\U}(X) : \ell^\infty (X)$ denotes the Banach space of all bounded
sequences in $X$, and $\N_{\U}(X)$ its subspace of sequences $(x_j)$ with
$\lim_{j\in\U} \|x_j\| =0$.
For any bounded sequence $(x_j)$ in $X$, we denote its equivalence class
in $X_{\U}$ by $[(x_j)]_{j=1}^\infty$.
For such an object $x$, we have $\|x\| = \lim_{j\in\U} \|x_j\|$.
Given $T\in \L(X,Y)$, the ultrapower  $T_{\U} :X_{\U} \to Y_{\U}$ is
defined by $T_{\U}(x) = [(Tx_j)]_{j=1}^\infty$ for all
$x= [(x_j)]_{j=1}^\infty$ in $X_{\U}$.
By the results in \cite{GM}, we may take as working definition:
$T$ is super weakly compact (resp. super Tauberian) provided $T_{\U}$ is
weakly compact (resp. Tauberian).
We note that (by results cited in
\cite{GM}), these definitions are independent of the chosen ultrafilter.

We first localize Proposition 2 and Corollary 4.

\proclaim{Proposition 19}
Let $T\in \L(X,Y)$.
Then the following are equivalent.
\roster
\item"1.)" $T$ is non-super weakly compact.
\item"2.)" There is a $\lambda\ge1$ so that for all $n$, there exist
$b_1,\ldots,b_n$ in $X$ with $(b_i)_{i=1}^n$ and $(Tb_i)_{i=1}^n$
$\lambda$-wide-$(s)$.
\item"3.)" There is a $\lambda\ge1$ so that for all $n$, there exist
$e_1,\ldots,e_n$ in $X$ with $(e_i)_{i=1}^n$ and $(Te_i)_{i=1}^n$
$\lambda$-wide-$(c)$.
\endroster
\endproclaim

\demo{Proof}
$2) \Leftrightarrow 3)$ follows by Proposition 10.
Now suppose first that condition 2) holds.
Fix $\U$ a non-trivial ultrafilter on $\nat$, and let $X_{\U},Y_{\U}$ be the
ultrapowers of $X$, $Y$ respectively.
For each $n$, let $b_1^n,\ldots,b_n^n$ be chosen in $X$ with
$(b_i^n)_{i=1}^n$ and $(Tb_i^n)_{i=1}^n$ $\lambda$-wide-$(s)$.
For $n<i$, set $b_j^n=0$, and
now define $(b_i)$ in $X_{\U}$ by $b_i= [(b_i^n)]_{n=1}^\infty$ for all $i$.
Then it follows that $(b_i)$ is $\lambda$-wide-$(s)$ and $(T_{\U}b_i)$
is also $\lambda$-wide-$(s)$.
But then $T_{\U}$ is not weakly compact by Proposition~2,
hence $T$ is not super weakly compact.
Conversely, if $T_{\U}$ is not weakly compact, then again applying
Proposition~2, there exists a sequence $(b_j)$ in $X_{\U}$ with
$(Tb_j)$ $\lambda$-wide-$(s)$ in $Y_{\U}$.
But now assuming $(b_j)$ and $(Tb_j)$ are both $\beta$-wide-$(s)$ and
letting $\lambda>\beta$, standard properties of ultraproducts allows us
to deduce the existence of finite sequences satisfying 2) for all $n$.\qed
\enddemo

\remark{Remark}
It follows moreover from Theorem $11'$ and the remark following that
in fact {\it $T$ is non-super weakly compact iff there is a $\lambda\ge1$ so
that for all $n$ and $\ep>0$ $T$ maps a normalized $(1+\ep)$-wide-$(s)$
sequence of length $n$ into a $\lambda$-wide-$(s)$ sequence\/}.
\endremark

We thus obtain immediately the following 

\proclaim{Corollary}
The following assertions are equivalent, for all given Banach spaces $X$.
\roster
\item"1)" $X$ is non-super reflexive.
\item"2)" There is a $\lambda\ge1$ so that $X$ contains $\lambda$-wide-$(s)$
sequences of arbitrarily large length.
\item"3)" For all $\ep >0$, $X$ contains normalized $(1+\ep)$-wide-$(s)$
sequences of arbitrarily large length.
\endroster
\endproclaim

\proclaim{Proposition 20}
Let $T\in \L(X,Y)$.
Then the following are equivalent.
\roster
\item"1.)" $T$ is non-super-Tauberian.
\item"2.)" There is a $\lambda\ge1$, so that for all $n$, all $\ep >0$,
there exist $b_1,\ldots,b_n$ in $X$ with $(b_i)_{i=1}^n$
$\lambda$-wide-$(s)$ and
$\diam \{Tb_i: 1\le i\le n\} <\ep $.
\item"3.)" For all $\ep >0$ and all $n$ there exist $b_1,\ldots,b_n$
in $X$ with $(b_i)_{i=1}^n$ $(1+\ep)$-wide-$(s)$ and
$\diam \{Tb_i: 1\le i\le n\} <\ep $.
\item"4.)" There is a $\lambda\ge1$  so that for all $n$, all $\ep >0$, there
exist $e_1,\ldots,e_n$ in $X$ with
$(e_i)_{i=1}^n$ $\lambda$-wide-$(c)$ and $\sum_{i=1}^n \|Te_i\|<\ep $.
\endroster
\endproclaim

\demo{Proof}
Again $2)\Leftrightarrow 4)$ follows by Proposition 10.

We now proceed with the same method as in the preceding case.
Again, suppose condition~2) holds.
Then for all $n$, we may choose $b_1^n,\ldots,b_n^n$ in $X$ with
$(b_i^n)_{i=1}^n$ $\lambda$-wide-$(s)$ and $\diam \{Tb_i^n: 1\le i\le n\}
< \frac1n$.
Again, letting $\U,X_{\U}, Y_{\U}$, and $T_{\U}$ be as in the preceding
proof, and again defining $b_i = [(b_i^n)]_{n=1}^\infty$, we have that
$(b_i)_{i=1}^\infty$ is a $\lambda$-wide-$(s)$ sequence in $X_{\U}$.
In this case, we have that $(T_{\U}b_i)$ is a constant sequence, hence
$T_{\U}$ is non Tauberian by Theorem~5.
(Thus $T$ is non-super-Tauberian by the results in \cite{GM}.)
Indeed, we have that for all $j$,
$$\|T_{\U} b_1 - T_{\U} b_j\| = \lim_{n\in \U} \|Tb_1^n -Tb_j^n\| =0\ .
\tag 33$$

Conversely, suppose $T$ is non-super-Tauberian.
Thus by \cite{GM}, $T_{\U}$ is non-Tauberian, so by Theorem~$5'$,
for all $\ep>0$, there is a $1+(\ep/2)$-wide-$(s)$ sequence $(b_j)$ in
$X_{\U}$ with $\diam \{Tb_i :1\le i\le n\}\le \ep/2$.
Again it follows from standard properties of ultraproducts that 3) holds.
Thus we have shown $2)\Rightarrow 1) \Rightarrow 3)$ and of course
$3)\Rightarrow 2)$ is trivial.\qed
\enddemo

\remark{Remarks}
1. A proof of the above result may also be obtained
by directly applying our Proposition~15 and Proposition~12 of \cite{GM}.

2. Our proof also yields (cf.\ the third remark following the proof of
Corollary~16 above) that
{\it if $T$ is non-super-Tauberian, then for all $\ep>0$ and $n$, there exists
$(b_1,\ldots,b_n)$ in $X$ which is $(2+\ep)$-wide-$(s)$ and
$\|Tb_i\| <\ep$ for all $1\le i\le n$.}
Actually, Proposition~12 of \cite{GM} and our Proposition~18 yield that
``2'' may be replaced by ``1'' in this assertion.
We may deduce this from
our work as follows:
Choose any $(b_1^n,\ldots,b_n^n)$ in $X$ with $(b_1^n,\ldots,b_n^n)$
3-wide-$(s)$ and $\|Tb_i^n\| < {1\over n}$ all $n$, and let
$b_j =[(b_j^n)]_{n=1}^\infty$ in $X_{\U}$ for all $j$
($\U$ a specified non-trivial ultrafilter on $N$).
But then $(b_j)$ is 3-wide-$(s)$ in $X_{\U}$, yet $T_{\U}b_j=0$ all $j$.
But this yields the result in \cite{GM} that $\ker T_{\U}$ is non-reflexive.
But then by Theorem~11,
we have that given $\ep>0$, there exists $(v_j)$ in $\ker T_{\U}$
with $(v_j)$ $(1+\ep)$-wide-$(s)$.
Of course now using standard ultraproduct techniques, we indeed obtain
that {\it if $T$ is non-super-Tauberian, then for all $n$ and $\ep$,
there exists $(b_i,\ldots,b_n)$ in $X$ $(1+\ep)$-wide-$(s)$ and
$\|Tb_i\|<\ep$ for all $i$}.
Since the $\ep$ here may be chosen arbitrarily, it also follows that in
fact $(b_1,\ldots,b_n)$ may be chosen normalized.
\endremark

Proposition 20, 3) also yields the result of D.G. Tacon \cite{T} that the
super-Tauberian operators form an open set in $\L(X,Y)$ (with the same
argument as the proof given in \cite{GM}).

\proclaim{Corollary}
The non-super-Tauberian operators from $X$ to $Y$ are a closed subspace
of $\L(X,Y)$.
\endproclaim

\demo{Proof}
Assume $T_n\to T$, $T_n$ non-super-Tauberian operators in  $\L(X,Y)$
for all $n$.
Now $\ep>0$ given, if $T_n$ is fixed with
$\|T_n-T\|<\ep$, for all $k$, choose $b_1,\ldots,b_k$
$(1+\ep)$-wide-$(s)$ with $\diam_{1\le i\le k} T_n(b_i) <\ep$.
But then $\diam_{1\le i\le k} T(b_i) < 2\ep (1+\ep) +\ep$.
Thus $T$ is non-super-Tauberian by Proposition~20 part~3).\qed
\enddemo

We now introduce a new class of operators, intermediate (as we shall see)
between the classes of Tauberian and super Tauberian operators.

\definition{Definition}
$T$ in $\L(X,Y)$ is called {\it strongly Tauberian\/}
provided its induced operator
$\widetilde T$ from $X^{**}/X$ to $Y^{**}/Y$ is an (into) isomorphism.
\enddefinition

\noindent
(In this definition, letting $\pi :X^{**}\to X^{**}/X$ denote the quotient
map, $\widetilde T$ is defined by $\widetilde T(\pi x^{**}) =
\pi (T^{**} x^{**})$.)
Evidently $T$ is Tauberian precisely when its induced operator
$\widetilde T$ is one--one.
The stronger property given by Definition~6 immediately yields that the
strongly Tauberian operators from an open semi-group; as we shall shortly
see, this class is also closed under such natural operations as taking
double adjoints.
We now list several equivalences for these operators.
(A word about notation: we regard $X\subset X^{**}$ and $X^{**}\subset
X^{4*}$; for $Y\subset X$, $Y^\bot = \{x^*\in X^*: x^* (y)= 0$ all $y\in Y\}$;
thus e.g., $X^\bot$ is a subspace of $X^{3*}$ so $X^{\bot\bot}$ is a
subspace of $X^{4*}$.
Recall that although $X^{\bot\bot}$ is canonically isometric to $X^{**}$,
$X^{\bot\bot}\cap X^{**} = X$.
Also, we denote $T^{***}$ by $T^{3*}$, $X^{***}$ by $X^{3*}$ etc.)

\proclaim{Theorem 21}
Let $T\in \L(X,Y)$.
Then the following are equivalent.
\roster
\item"1)" $T$ is strongly Tauberian.
\item"2)" There is a $\delta >0$ so that for all $x^{**}\in X^{**}$,
$\dist (T^{**}x^{**},Y) \ge \delta\dist (x^{**},X)$.
\item"3)" $T^{3*} Y^\bot = X^\bot$
\item"4)" $T^{4*} | X^{*\bot}$ is an isomorphism.
\item"5)" $T^{**}$ is strongly Tauberian.
\item"6)" For all $\lambda>1$, there exists a $\beta>1$ so that every
(infinite) $\lambda$-wide-$(s)$ sequence in $X$ has a subsequence which
is mapped into a $\beta$-wide-$(s)$ sequence by $T$.
\item"7)" Same as 6), except replace ``for all $\lambda>1$'') by
``there exists $\lambda>1$''.
\endroster
\endproclaim

\remark{Remark}
It follows immediately from the definition that if $T,S\in \L(X,Y)$
are given with $T$ strongly Tauberian and $S$ weakly compact, then
$T+S$ is strongly Tauberian, since $\widetilde{T+S} = \widetilde T+
\tilde S= \widetilde T$.
We obtain a simultaneous generalization of this observation, together
with the openness of the set of strongly Tauberian operators in $\L(X,Y)$, 
as follows:
{\it Given $T$ strongly Tauberian, there exists $\ep>0$ so that if
$S\in \L(X,Y)$ satisfies $\|\tilde S\|<\ep$, then $S+T$ is strongly
Tauberian.}
Indeed, simply choose $\delta$ as in Theorem~21 part 2), then any $\ep
<\delta$ works.
\endremark

\demo{Proof of Theorem 21}
$1)\Leftrightarrow 2)$ Immediate from Definition~6.

$1)\Leftrightarrow 3)$ $(X^{**}/X)^*$ is canonically identified with
$X^\bot$; now we always have that $T^{3*}Y^\bot \subset X^\bot$ and
moreover $(\widetilde T)^*$ may be canonically identified with
$T^{3*}|Y^\bot$.
Now this equivalence follows from the observation that an operator is an
isomorphism iff its adjoint is surjective.
\enddemo

For the next two equivalences, it is convenient to isolate out the
following elementary fact.

\proclaim{Lemma}
Let $Z,W$, be Banach spaces, $Z_i$, $W_i$ be closed linear subspaces of
$Z,W$ respectively, $1\le i\le 2$, $S:Z\to W$ be a given bounded linear
operator with $SZ_i\subset W_i$, $1\le i\le 2$, and $Z= Z_1\oplus Z_2$,
$W= W_1\oplus W_2$.
Let $\tilde S :Z/Z_1\to W/W_1$ be the induced operator defined by
$\tilde S(\pi z)= \pi (Sz)$, $\pi$ the appropriate quotient map.
Then $\tilde S$ is an isomorphism iff $S|Z_2$ is an isomorphism.
\endproclaim

$1)\Leftrightarrow 4)$ 
We apply the Lemma to $Z$ (resp. $W$) $=X^{4*}$ (resp. $Y^{4*}$),
$Z_1$ (resp. $W_1$) $= X^{\bot\bot}$ (resp. $Y^{\bot\bot}$),
$Z_2 = X^{*\bot}$ (resp. $Y^{*\bot}$), and $S= T^{4*}$.
Now we have, by standard Banach space theory, that $X^{3*} = X^\bot \oplus
X^*$ hence $X^{4*} = X^{\bot\bot} \oplus X^{*\bot}$.
Thus the hypotheses of the Lemma are fulfilled.
Of course for any operator $U$, $U$ is an isomorphism iff $U^{**}$ is an
isomorphism.
Now first assuming $T$ is strongly Tauberian, we thus have that $(\widetilde
T)^{**}$ is an isomorphism.
But $(X^{**}/X)^{**}$ may be canonically identified with $X^{4*}/X^{\bot\bot}$
and $(\widetilde T)^{**}$ with $\tilde S$.
Hence by the Lemma, $\tilde S|Z_2 = T^{4*}|X^{*\bot}$ is an isomorphism.
But conversely, $T^{4*}|X^{*\bot}$ an isomorphism implies
$(\widetilde T)^{**}$ is an isomorphism, again by the Lemma, so $T$ is
strongly Tauberian.

$4)\Leftrightarrow 5)$
Now we apply our Lemma to the fact that $X^{4*}= (X^*)^{3*}  = X^{*\bot}
\oplus X^{**}$.
Thus we obtain that $T^{4*}|X^{*\bot}$ is an isomorphism iff
$\widetilde{T^{**}}$ is an isomorphism.

$2) \Rightarrow 6)$
Let $T$ be strongly Tauberian, and choose $\delta >0$ satisfying 2).
Now let $(x_j)$ be a $\lambda$-wide-$(s)$ sequence in $X$, and let $G$
be a $w^*$-cluster point of $(x_j)$ in $X^{**}$.
By an argument in the proofs of Corollaries~16, 17, we have that
$$\dist (G,X) \ge {1\over 2\lambda}\ .
\tag 34$$
Then $T^{**}G$ is a $w^*$-cluster point of $(Tx_j)$, so by 2),
$$\dist (T^{**}G,Y)\ge {\delta\over 2\lambda}\ .
\tag 35$$
Hence by the Corollary to Theorem 12, $(Tx_j)$ has a $\beta$-wide-$(s)$
subsequence, where ($\ep>0$ given)
$$\beta = \left({2\lambda\over\delta}
\vee {\lambda^2 +{\delta\over2}\over\delta}
\vee \lambda\right) +\ep\ .
\tag 36$$

$6)\Rightarrow 7)$ is trivial, so it remains to prove $7)\Rightarrow 1)$.
This is an immediate consequence of the following equivalence.

\proclaim{Proposition 22}
$T\in \L(X,Y)$ is non-strongly Tauberian iff for all $\ep>0$, there
exists $(x_j)$ a $(1+\ep)$-wide-$(s)$ sequence in $X$ with
$\diam \{Tx_1,Tx_2,\ldots\} <\ep$.
\endproclaim

\demo{Proof of Proposition 22}
It suffices to prove the direct implication.
Indeed, if the second assertion of Proposition 22 holds, then (since
trivially any $1+\eta$-wide-$(s)$ sequence is $1+\ep$-wide-$(s)$ if
$\eta<\ep$), 6) of Theorem~21 fails, whence by $1)\Rightarrow 6)$ of the
latter, $T$ is non-strongly Tauberian.

Now assume $T$ is non-strongly Tauberian, and let $0<\eta <1$.
Then we may choose $G$ in $X^{**}$ satisfying
$$\|G\| <1\ ,\ \dist (G,X) >1-\eta\ ,\ \text{ and }\ \dist (T^{**}G,Y)
< \eta \ .
\tag 37$$
Of course then we may also choose $y\in X$ with
$$\|T^{**} G-y\| <\eta\ .
\tag 38$$
\enddemo

Now standard techniques (cf.\ \cite{R1}) yield

\proclaim{Lemma 23}
Let $L= \{x\in X: \|x\| <1$ and $\|Tx-y\|<\eta\}$.
Then $G\in \tilde L$.
\endproclaim

(For $M\subset B$ a Banach space, $\widetilde M$ denotes the
$w^*$-closure of $M$ in $B^{**}$.)

\demo{Proof of Lemma 23}
If not, by the Hahn-Banach separation theorem, choose $x^*\in X^*$ and
$a<b \eqdf G(x^*)$ so that
$$x^*(\ell) \le a\ \text{ for all }\ \ \ell\in L\ .
\tag 39$$
Now it follows that setting $W= \{x\in X:\|x\|<1$ and $x^*(x)> {a+b\over2}\}$, 
then $G\in \widetilde W$ and of course $W\cap L=\emptyset$.
Thus $W$ is a convex set so that
$$T^{**}G\in \widetilde{TW}\ \text{ and }\ \|Tw-y\| \ge \eta\ \text{ all }\
w\in W\ .
\tag 40$$
Then again by the  Hahn-Banach theorem, there exists a $y^*$ in $Y^*$
with $\|y^*\| =1$ and
$$y^* (Tw-y) \ge \eta\ \text{ all }\ w\in W\ .
\tag 41$$
But since $T^{**}G-y$ is in the $w^*$-closure of $TW-y$, we obtain that
$$\|T^{**}G-y\| \ge \langle y^*,T^{**}G-y \rangle \ge \eta\ ,
\tag 42$$
contradicting (38). (Here we have assumed $L\ne\emptyset$.
However if $L=\emptyset$, instead let $W= \{ x\in X:\|x\|<1\}$; 
now (40) holds, and the rest of the argument following (40) again yields
a contradiction to (38).\qed
\enddemo

We may now complete the proof of Proposition 22 as follows.
Let $\ep >0$, and let $\eta >0$ be chosen, with ${1\over 1-\eta}<\ep $.
First, the proof of Theorem~12 yields that we may choose $(b_j)$ a
${1\over 1-\eta}$-triangular sequence in $L$.
(The proof of the existence of $(b_j)$ in $L$, $(f_j)$ in $X^*$
satisfying (13)--(16) does not require the separability of $X$.)
Here, we are just using that $G\in \tilde L$ and $\dist (G,X) > 1-\eta$.
Then by Corollary~(16), $(b_j)$ has a
$[{({1\over1-\eta})^2 +1\over 2} +\eta]$-wide-$(s)$ subsequence $(x_j)$.
Of course we are now finished, as long as
$${({1\over 1-\eta})^2 +1\over 2} +\eta < 1+\ep \ \text{ and }\
\eta < {\ep\over2}\ .$$
Indeed, then since $x_j\in L$ for all $j$, $\|Tx_i-Tx_j\| < 2\eta <\ep$
all $i,j$.\qed

\remark{Remark}
The proof of $2)\Rightarrow 6)$ yields a non-linear estimate for the
dependence of $\beta$ on $\lambda$  (36).
However using Theorem~12 itself, we obtain that
{\it if $\delta$ is as given in $2)$ of Theorem~21, then given $\ep>0$,
every infinite $\lambda$-triangular sequence in $X$ has a subsequence
which is mapped into a $({2\lambda\over\delta}+\ep)$-triangular
sequence by\/} $T$ which is also $\beta$-wide-$(s)$, with $\beta$ as in (36).
Thus working with $\lambda$-triangular sequences, we recapture the best
``$\delta$'' in 2), to within a factor of 2.
\endremark

Of course Proposition 22 and Proposition 20, part 3) have the following
immediate consequence.

\proclaim{Corollary 24}
Every super-Tauberian operator is strongly Tauberian.
\endproclaim

\remark{Remark} 
We thus have for a given $T\in \L(X,Y)$,
$$T\text{ super-Tauberian } \Rightarrow T \text{ strongly Tauberian }
\Rightarrow T \text{ Tauberian.}$$
It is easily seen these implications are strict (in general).
Indeed, there exist reflexive Banach spaces $X$ which admit non-super
Tauberian operators on them (e.g., $X= (\bigoplus_{n=1}^\infty \ell_n^1)_2$; 
of course any operator on a reflexive Banach space is
strongly Tauberian.
If $X=(\bigoplus_{n=1}^\infty c_0)_2$, and $T((x_n)) = ({1\over n}x_n)$
all $(x_n)\in X$, then $T$ is Tauberian and
any open neighborhood of $T$ in $\L(X)$ contains a non-Tauberian
operator, namely $Sx_n= \frac1n x_n$ for $n\le N$, $S(x_n)=0$ for $n>N$,
for suitable $N$.
Hence $T$ is non-strongly Tauberian since the strongly Tauberian operators
are an open set in $\L(X)$, contained in the Tauberian ones.
(I am indebted to A.~Mart{\'\i}nez for this  conceptual proof that $T$ is
non-strongly Tauberian.)
It is also known there exist Tauberian operators $T$ with $T^{**}$
non-Tauberian \cite{AG}.
\endremark

Our next result, localizing Corollary 6, yields a ``direct'' characterization
of super-Tauberian operators.

\proclaim{Corollary 25}
Let $T\in \L(X,Y)$.
Then the following are equivalent.
\roster
\item"1)" $T$ is super-Tauberian
\item"2)" For all $\lambda\ge1$ there exists $\beta\ge1$ so that for all
positive integers $k$, there exists $n>k$ so that if $(x_1,\ldots,x_n)$
is a $\lambda$-wide-$(s)$ sequence in $X$; then $(Tx_{j_1},\ldots,Tx_{j_k})$ is
$\beta$-wide-$(s)$ for some $j_1<j_2<\cdots j_k \le n$.
\item"3)" Same as {\rm 2)}, except replace ``For all $\lambda \ge1$''
``there exists $\lambda\ge1$''.
\endroster
\endproclaim

\remark{Remark}
This corollary immediately yields ``directly'' the semi-group property:
$ST$ is super Tauberian if $T$, $S$ are.
\endremark

Now fix $\lambda\ge1$, and let $\U$, $X_{\U}$, $Y_{\U}$, and $T_{\U}$
be as in the proof of Proposition~9.
Since $T$ is super-Tauberian, it follows easily that also $T_{\U}$ is
super-Tauberian.
Hence by Corollary~24, $T_{\U}$ is strongly Tauberian.
Thus we may choose $\beta>1$ so that
$$\text{Theorem 21  holds for ``$X$'' $=X_{\U}$, ``$T$'' $= T_{\U}$.}
\tag 43$$
Now we claim that $2\beta$ works for the Corollary.
If not, then there exists a positive $k$ so that for all $n$, we may
choose $(x_1^n,\ldots,x_n^n)$ $\lambda$-wide-$(s)$ in $X$, so that no
subsequence of $(Tx_1^n,\ldots,Tx_n^n)$ of length $k$ is $2\beta$-wide-$(s)$.
Now as usual, define $(x_j)$ in $X_{\U}$ by $x_j = [(x_j^n)]_{n=1}^\infty$
for all $j$.
Then it follows that $(x_j)$ is $\lambda$-wide-$(s)$ in $X_{\U}$.
Hence by (43), there exists $j_1<j_2<\cdots$ with $(T_{\U}x_{j_i})$ a
$\beta$-wide-$(s)$ sequence.
But then for $n$ sufficiently large, $(Tx_{j_1}^n,\ldots, Tx_{j_k}^n)$ is
$2\beta$-wide-$(s)$, a contradiction.\qed

We conclude with a perturbation result whose proof follows from
ultraproduct considerations, which generalizes Proposition~18.
This result is a localization of Corollary~17 (cf.\ the remark following
its proof).
We say that for sequences $(x_j)$ and $(y_j)$ in a Banach space,
$(y_j)$ is an $\ep$-perturbation of $(x_j)$ if $\|x_j-y_j\| \le \ep$ all $j$.

\proclaim{Proposition 26}
Let $\lambda\ge1$ be given.
There exist $\beta\ge1$ and $\ep>0$ so that for all $k$, there is an
$n>k$ so that every $\ep$-perturbation of a $\lambda$-triangular sequence
of length $n$ has a $\beta$-wide-$(s)$ subsequence of length $k$.
\endproclaim

\demo{Proof}
Let $\beta$ and $\ep$ be chosen satisfying the conclusion of Corollary~17,
and let $\beta'>\beta$.
We shall show that $\beta'$ and $\ep$ satisfy the conclusion of
Proposition~26.
Were this false, by simply taking the $c_0$-sum of some finite-dimensional
Banach spaces, we may choose a Banach space $X$ and for every $n$,
sequences $(b_j^n)_{j=1}^n$ and $(p_j^n)_{j=1}^n$ in $X$ with
$(b_j^n)_{j=1}^n$ $\lambda$-triangular, $\| p_j^n\| \le \ep$ all $j$,
yet $(b_j^n +p_j^n)$ has no $\beta'$-wide-$(s)$ subsequence of length~$k$.
But then letting $\U$ and $X_{\U}$ be as before,
and defining $b_j=[(b_j^n)]_{n=1}^\infty$,
$p_j = [(p_j^n)]_{n=1}^\infty$ for all $j$, then $(b_j)$ is
$\lambda$-triangular in $X_{\U}$ and $\|p_j\|\le \ep$ for all $j$,
hence $(b_j+p_j)$ has a $\beta$-wide-$(s)$ subsequence in $X_{\U}$, by
Corollary~17.
But then for $n$ sufficiently large, $(b_j^n +p_j^n)_{j=1}^n$ has a
$\beta'$-wide-$(s)$ subsequence of length $k$, contradicting our
assumption.~\qed
\enddemo

\enddocument